\documentclass[12pt,a4paper]{article}
\usepackage[a4paper,vmargin={20mm,20mm},hmargin={27mm,27mm}]{geometry}
\usepackage[ruled,vlined]{algorithm2e}
\usepackage{dirtytalk}
\usepackage{comment} 
\usepackage{caption}
\usepackage{amsmath}
\usepackage{float}
\usepackage[utf8x]{inputenc}
\usepackage{amsmath,amssymb}
\usepackage{xcolor}
\usepackage{mathtools}
\usepackage{tabulary}
\usepackage{amsthm}
\usepackage{latexsym}
\usepackage[mathcal]{eucal}
\usepackage{amscd}
\usepackage{graphicx}
\usepackage{amsfonts}
\usepackage{amscd}
\usepackage{MnSymbol}
\usepackage{setspace}
\newtheorem{definition}{Definition}[section]
\usepackage[english]{babel}
\usepackage{listings}
\usepackage{calligra}
\usepackage{titlesec}

\newtheorem{theorem}{Theorem}[section]
\newtheorem{example}[theorem]{Example}
\newtheorem{exercise}[theorem]{Exercise}
\newtheorem{lemma}[theorem]{Lemma}

\newtheorem{proposition}[theorem]{Proposition}
\newtheorem{remark}[theorem]{Remark}

\newcommand{\bt}{\begin{theorem}}
\newcommand{\et}{\end{theorem}}
\newcommand{\bl}{\begin{lemma}}
\newcommand{\el}{\end{lemma}}
\newcommand{\bexc}{\begin{exercise}}
\newcommand{\eexc}{\end{exercise}}
\newcommand{\bpr}{\begin{proposition}}
\newcommand{\epr}{\end{proposition}}
\newcommand{\bex}{\begin{example}}
\newcommand{\eex}{\end{example}}
\newcommand{\bc}{\begin{corollary}}
\newcommand{\ec}{\end{corollary}}
\newcommand{\bo}{\begin{proof}}
\newcommand{\eo}{\end{proof}}
\newcommand{\bd}{\begin{definition}}
\newcommand{\ed}{\end{definition}}
\newcommand{\br}{\begin{remark}}
\newcommand{\er}{\end{remark}}
\newcommand{\be}{\begin{enumerate}}
\newcommand{\ee}{\end{enumerate}}

\newcommand{\K}{\mathcal{K}}

\newcommand{\Z}{\mathbb{Z}}

\begin{document}

\title{
    {\textbf{Coding billiards in hyperbolic 3-space}}\\ \vspace{0.1cm}
    \author{\large Pradeep Singh \vspace{0.1cm}\\
    \normalsize Department of Mathematics \vspace{-0.15cm}\\ \normalsize  Indian Institute of Technology Delhi \vspace{-0.15cm}\\ \normalsize pradeep.singh@maths.iitd.ac.in}
}

\date{\today}

\maketitle

\begin{abstract}
In this paper, we extend the scope of symbolic dynamics to encompass a specific class of ideal polyhedrons in the 3-dimensional hyperbolic space, marking an important step forward in the exploration of dynamical systems in non-Euclidean spaces. Within the context of billiard dynamics, we construct a novel coding system for these ideal polyhedrons, thereby discretizing their state and time space into symbolic representations. This  paper distinguishes itself through the establishment of a conjugacy between the space of pointed billiard trajectories and the associated shift space of codes. A crucial finding herein is the observation that the closure of the related shift space emerges as a subshift of finite type (SFT), elucidating the structural aspects and asymptotic behaviour of these systems. 
\end{abstract}

\bigskip


\section{Introduction}

The realm of Symbolic Dynamics constitutes a significant subset within the expansive theory of dynamical systems. This field focuses on delineating a specific dynamical system by discretizing both its state and time space. The discretized regions within the state space are assigned labels, thereby enabling the inherent dynamics to generate corresponding sequences of symbols. If an equivalence between this symbolic space and the original dynamical system can be established, the system can be marked with a simple computational structure through the symbolic space. This facilitates the application of a multitude of analytic tools designed for symbolic spaces to the dynamical system, enhancing our understanding of the latter.\\

Billiard dynamics are frequently studied across physics and mathematics, primarily in Euclidean spaces. However, the implications of gravity naturally lead us to explore non-Euclidean spaces. Practically, billiard theory often models the price action of various financial instruments, with physical and notional constraints acting as the boundaries of the billiard table. Here our central interest lies in understanding how various billiard trajectories relate to each other asymptotically for a given billiard table configuration. As a result, the dynamics of sets take precedence over individual dynamics. Establishing a link between the original dynamical system and its hyperspace (i.e., the space of subsets) can greatly enhance our understanding of these dynamics. Numerous studies have explored this connection, outlining the correlation between individual and set-valued dynamics \cite{banks, sharma1, sharma2}. In this article, we will extend this understanding to billiard systems in 3-dimensional hyperbolic space, leveraging the established theory of pointed geodesics \cite{nagar}. Further, we will provide a symbolic encoding of such systems.

The problem of coding billiards in the hyperbolic plane has been extensively studied \cite{castle,ullmo,ullmo2, nagar}. While the 2-dimensional problem is supplemented by the geometric intuition one gains from the two-dimensional models of the hyperbolic plane, the case of higher dimensions lacks similar intuitive models. As a result, work in 3-dimensions and higher has been primarily focused on the issue involving relativistic billiards \cite{deryabin1,deryabin2,deryabin3,pustylnikov}.\\

This article seeks to investigate the coding of billiards in a class of polyhedra in 3-dimensional hyperbolic space. Our polyhedra are situated in the Poincaré ball model $\mathbb{B}^3= \{(x,y,z) \in \mathbb{R}^3 : x^2+y^2+z^2 <1 \}$ with vertices located on the boundary $\partial \mathbb{B}^3$. Consequently, all vertices of such polyhedra sit at infinity under hyperbolic distance. We denote this class as \say{ideal.} For a given ideal polyhedron, we define a billiard trajectory as a point inside the polyhedron moving with uniform speed along geodesics until it hits a face of the polyhedron. It then reflects elastically according to the laws of specular reflection in the hyperbolic plane, as defined by the incident geodesic and the normal geodesic at the point of contact. After the collision, it follows the reflected geodesic until it encounters another face of the polyhedron, at which point the process repeats indefinitely into the future. The same can be described in the past by considering the point moving in the opposite direction. This results in a bi-infinite collection of geodesic segments, which we term a \say{billiard trajectory.} It is crucial to note that we exclude points that hit an edge or a vertex of the polyhedron either in the past or future. Our main interest lies in non-truncated pasts and futures. We arbitrarily label the faces of the polyhedron. Upon isolating a billiard trajectory, it naturally generates a bi-infinite sequence of symbols corresponding to the faces hit in sequence. We select a geodesic segment from the billiard trajectory, termed the \say{base} of the \say{pointed billiard trajectory.} This leads to the determination of a \say{base symbol} in the corresponding bi-sequence of symbols. We term this new concept as the \say{code for the pointed billiard trajectory.} The primary objective of this article is to classify the pointed billiard trajectories for the class of ideal polyhedra via certain grammar rules on the associated shift space.\\

In Section 2, we introduce the fundamental concepts necessary to address the problem at hand, as well as describe the billiard map pertinent to billiards in hyperbolic space. Following this, Section 3.1 delves into the concept of pointed geodesics for billiards within ideal polyhedra. We define a metric on the space of pointed geodesics that bears topological equivalence to the Hausdorff metric. Then, in Section 3.2, we present the core result of this work, which sets forth the coding rules for billiard trajectories within a class of ideal polyhedra in hyperbolic space. Using these coding rules, we demonstrate the conjugacy between the space of pointed geodesics and the associated symbolic space.

\section{Preliminaries}
In this section we lay down some basic notions to be used further. The hyperbolic 3-space, denoted as $\mathbb{H}^3$, is a 3-dimensional Riemannian manifold that carries a constant negative curvature. 
The notion of curvature in a Riemannian manifold is central to the definition of a hyperbolic space and characterizes the amount by which the geometry of the space deviates from that of the flat Euclidean space. $\mathbb{H}^3$ carries forward the properties of a hyperbolic plane, but in three dimensions, thereby offering three degrees of freedom. The  hyperbolic plane cannot be embedded in $\mathbb{R}^2$ and thereby is studied using various models, most common being the Poincar\'e half plane model denoted $\mathbb{H}^2$ and the Poincar\'e     disc model denoted $\mathbb{D}^2$, see e.g., \cite{anderson}. Same holds true for  the hyperbolic 3-space which is dealt with two commonly used models - Poincar\'e half space model denoted $\mathbb{H}^3$ and Poincar\'e ball model denoted $\mathbb{B}^3$.\\

The Poincar\'e half space model is defined with the underlying space as $\mathbb{H}^3 = \{(x,y,z) \in \mathbb{R}^3 : z >0\}$ and an attached metric given by 
\begin{equation}
    ds^2 = \dfrac{dx^2 +dy^2+dz^2}{z^2},
    \end{equation}
    which is the Euclidean metric adjusted to ensure constant negative curvature. This is consistent with the general definition of a Riemannian metric, which is a type of smoothly varying inner product on the tangent space of a manifold.
We think of the boundary of $\mathbb{H}^3$ as the complex plane $\mathbb{C}$ with the natural embedding and a point at infinity. This allows us to extract the isometry group of $\mathbb{H}^3$
in terms of the available isometry group of $\mathbb{C}$. If $\gamma : [a,b] \rightarrow \mathbb{H}^3$ is a path in $\mathbb{H}^3$ given by $\gamma(t) = (x(t),y(t),z(t))$ with $t \in [a,b]$, then
\begin{equation}
l_{\mathbb{H}^3}(\gamma) = \int_{\gamma} \dfrac{\sqrt{dx^2+dy^2+dz^2}}{z}
\end{equation}
gives its length. 
The distance between two points $A$ and $B$ is given by  $d_{\mathbb{H}^3}(A,B)= \inf(l_{\mathbb{H}^3}(\gamma) )$ with the infimum running over all the paths originating from $A$ and terminating at $B$.\\

The underlying space for Poincar\'e ball model is $\mathbb{B}^3 = \{(x,y,z) \in \mathbb{R}^3 : x^2 +y^2+ z^2 < 1 \}$. Its metric 
\begin{equation}
ds^2 = \frac{4(dx^2+dy^2+dz^2)}{(1-x^2-y^2-z^2)^2}
\end{equation}
is designed such that the geodesics are the chords of the unit ball, providing an intuitive analog to Euclidean geometry. If $\gamma : [a,b] \rightarrow \mathbb{B}^3$ is a path in $\mathbb{B}^3$ that is given by $\gamma(t) = (x(t),y(t),z(t))$, then its length is given  by 
\begin{equation}
l_{\mathbb{B}^3}(\gamma) = \int_{\gamma} \dfrac{2 \sqrt{dx^2+dy^2+dz^2}}{1-(x^2+y^2+z^2)}.
\end{equation}
With points $A,B \in \mathbb{B}^3$, their distance is $d_{\mathbb{B}^3}(A,B)= \inf(l_{\mathbb{B}^3}(\gamma) )$ with the infimum running over all the paths originating from $A$ and terminating at $B$.
We can establish isometries between $\mathbb{H}^3$ and $\mathbb{B}^3$ analogous to the Cayley map between $\mathbb{H}^2 $ and $\mathbb{D}^2$. For a more detailed account on the hyperbolic space,  \cite{bonahon, marden} are suggested and for an intuitive understanding of the same, \cite{thurston, thurston2} are recommended. \\

Geodesics are the shortest paths between two points in a Riemannian manifold, and serve as the hyperbolic equivalent of \say{straight lines} in Euclidean space.  Roughly, these are the locally distance minimising curves of the space. For $\mathbb{H}^3$, the geodesics are the euclidean straight lines perpendicular to the $\{z=0\}$ plane and the euclidean semicircles hitting the $\{z=0\}$ plane orthogonally. The boundary of $\mathbb{H}^3$ denoted $\partial{\mathbb{H}^3}$  comprises of the plane $\{z=0\}$ and a point at infinity. The objects analogous to planes in Euclidean geometry in case of Hyperbolic geometry are called \emph{hyperbolic planes}. A hyperbolic plane in $\mathbb{H}^3$ model is either a euclidean plane hitting perpendicularly at the $\{z=0\}$ plane  or a euclidean hemisphere with its center lying on $\{z=0\}$ plane.

For $\mathbb{B}^3$, the geodesics are the euclidean straight lines passing through origin and the portions of  euclidean circles lying in $\mathbb{B}^3$ and hitting    $\partial{\mathbb{B}^3}$ orthogonally.
 The boundary of $\mathbb{B}^3$ comprises of the euclidean unit sphere 
centered at the origin. A hyperbolic plane in $\mathbb{B}^3$ model is either part of a euclidean plane passing through the origin  lying in $\mathbb{B}^3$  or part of a euclidean hemisphere inside $\mathbb{B}^3$ with its center lying on $\partial{\mathbb{B}^3}$.\\

A  \emph{convex} subset \emph{A} of the hyperbolic space is the one for which the closed line segment $l_{xy}$, joining $x$ to $y$ is contained in \emph{A} for each pair of distinct points $x$ and $y$ in \emph{A}. All hyperbolic lines, hyperbolic rays, hyperbolic line segments, hyperbolic planes, connected subsets of hyperbolic planes are convex.
Consider a hyperbolic plane $P$, the complement of $P$ in the hyperbolic space has two connected components, which are the two \emph{open half spheres determined by P}. We refer to $P$ as the \emph{boundary plane} for the half spheres it determines. Open half spheres and their closures are also convex in hyperbolic space.\\

Let $S = \{S_{\alpha}\}_{\alpha \in \Lambda}$ be a collection of half spaces in the hyperbolic space and for each $\alpha \in \Lambda$, let $P_{\alpha}$ be the bounding plane for $S_{\alpha}$. The collection $S$ is \emph{locally finite}
 if for each point $a$ in the hyperbolic space, $\exists \ \epsilon >0$ such that only finitely many bounding planes $P_{\alpha}$ intersect the open hyperbolic ball  $U_{\epsilon}(a)$. A \emph{hyperbolic polyhedron} is a closed convex set in the hyperbolic space that can be represented as the intersection of a locally finite collection of closed half spaces. With this definition some \emph{degenerate}  objects still enter the picture. Simplest example that we can take is that of any hyperbolic plane $P$ which satisfies the criteria but has empty interior. To disallow this we must ensure that a hyperbolic polyhedron is \emph{nondegenerate} which are the ones with nonempty interior. Here we will consider only nondegenerate hyperbolic polyhedrons. Suppose $\Pi$ be a hyperbolic polyhedron and $P$ a hyperbolic plane such that $\Pi$ intersects $P$ and $\Pi$ is contained in a closed half space determined by $P$. If the intersection $\Pi \cap P$ ends up in a point, we call it a $vertex$ of $\Pi$. If $\Pi \cap P$ is a hyperbolic line or a hyperbolic line segment, it is called an $edge$ of $\Pi$. If $\Pi \cap P$ is $P$ itself or its closed 2-dimensional connected subset, then it is called a $face$ of $\Pi$.\\

 Let $\Pi$ be a hyperbolic polyhedron and $e$ be an edge of $\Pi$  that is the intersection of two faces $F_1$ and $F_2$ of $\Pi$. Let $P_i$ be the hyperbolic plane containing $F_i\ \forall \ i =1,2$. Then $P_1 \cup P_2$ divides the hyperbolic space into four connected components, out of which one contains $\Pi$. The \emph{interior angle} of $\Pi$ at $e$ is the angle between $P_1$ and $P_2$ measured in the component of $P_1 \cup P_2$ containing $\Pi$. $\Pi$ has an \emph{ideal vertex} at a vertex $v$ if there are two adjacent faces of $\Pi$ that share $v$ as an endpoint on the boundary of the hyperbolic space. A finite-faced polyhedron $\Pi$ in the hyperbolic space is called $reasonable$ if it does not contain an open half-space. A \emph{hyperbolic k-hedron} is a reasonable hyperbolic polyhedron with $k$ faces.
A \emph{compact polyhedron} is a hyperbolic polyhedron whose all vertices are in the hyperbolic space. For $k \geq 4$, an \emph{ideal k-hedron} is a reasonable hyperbolic polyhedron $\Pi$ that has $k$ faces with all the vertices lying on the boundary of the hyperbolic space. A polyhedron is $isogonal$ if all its vertices are of same type i.e., each vertex is surrounded by the same number of faces and has same solid angle. Some features of an ideal polyhedron call for attention. The primary one being that the faces of an ideal polyhedron are themselves ideal polygons on which the billiard codes have been studied in \cite{castle, ullmo, nagar}. An $k$-faced polyhedron has $(k+4)/2$ vertices whose surface can be divided into $k$ faces consisting of ideal triangles. In particular, ideal polyhedrons come only as even-faced. Since area of each ideal triangle is $\pi$, the total surface area of the polyhedron itself comes out to be $k \pi$.  

\bigskip

\emph{This allows one to pose an interesting question of the relationship between the space of codes corresponding to the billiards on faces of a polyhedron with the space of codes for the bodily billiards that we seek to study in this article}.   

\bigskip
Moreover, not all euclidean polyhedrons have their ideal hyperbolic counterparts. A necessary condition for the existence of an ideal polyhedra corresponding to a euclidean polyhedra is if its vertices can be simultaneously put on a circumscribing sphere. If two ideal polyhedrons have the same number of vertices then they have the same surface area.         

 \begin{figure}[ht!]
\centering
\includegraphics[width=14.6cm,height=14.6cm]
{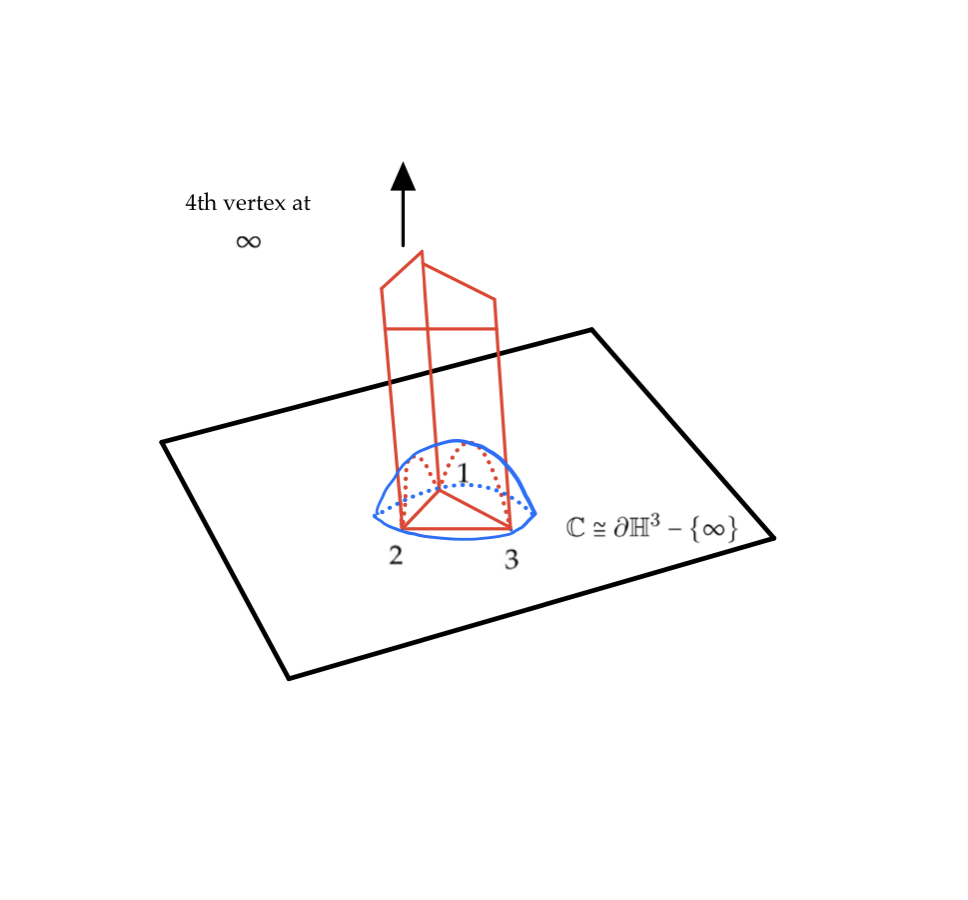}
\caption{An ideal polyhedron in $\mathbb{H}$}
\end{figure}

All plane face angles and solid angles at the vertices of an ideal polyhedron are zero. The angles between the faces (dihedral angles) are nonzero with their supplementary counterparts summing up to $2 \pi$. Dihedral angles are the angles formed by two intersecting hyperplanes (the higher-dimensional analog of planes in 3-space) in $\mathbb{H}^3$. They are used in studying polyhedra in hyperbolic 3-space because they can be used to measure the \say{bend} between faces of a polyhedron. Thus for an ideal isogonal tetrahedron, the dihedral angles come out to be $\pi/3$.  In this article, we are primarily focused on ideal isogonal polyhedrons. Figure 1 depicts a typical ideal tetrahedron. For further details on low-dimensional geometry, \cite{bonahon,marden, loring, thurston} are recommended.  A classical discussion on the related theory of fields can be found in \cite{landau}.

\subsection{Unfolding the Billiard Table}

In the following, we consider the Poincaré ball model $\mathbb{B}^3$ of hyperbolic space and focus our attention on the class of ideal hyperbolic polyhedrons, which we will henceforth refer to simply as \emph{ideal polyhedrons}.\\

Let us consider an ideal polyhedron $\Pi$ in $\mathbb{B}^3$ with $k$ faces. In the context of hyperbolic billiards, a billiard trajectory within this polyhedron $\Pi$ is defined as a directed geodesic flow that connects each pair of specular reflections off the faces (excluding the vertices and edges) of $\Pi$.\\

To facilitate our discussion, we parameterize the boundary of the Poincaré ball $\mathbb{B}^3$ using the azimuthal angle $\theta$ and the polar angle $\phi$. With this parameterization, the boundary of $\mathbb{B}^3$ can be considered as a subset of $\mathbb{R}^3$.

A directed geodesic can be represented as a pair $((\theta_i , \phi_i),(\theta_f, \phi_f) )$, where $(\theta_i , \phi_i)$ and $(\theta_f , \phi_f)$ are the intercepts made by the geodesic on $\partial \mathbb{B}^3$, with the direction originating from $(\theta_i , \phi_i) $ and leading to $(\theta_f , \phi_f) $. With these definitions, we equip $\partial \mathbb{B}^3$ with a metric that records the \emph{great-circle distance} between any two points on $\partial \mathbb{B}^3$ in terms of their latitudes and longitudes.\\

Next, we introduce the concept of the bounce map, which describes the evolution of geodesic arcs upon reflection off the faces of the polyhedron. Consider a pair $((\theta_i , \phi_i), (\theta_f , \phi_f)) $ that remains constant between any two consecutive reflections from the faces. Upon reflection, this pair transforms to a new pair $((\theta'_i , \phi'_i), (\theta'_f , \phi'_f)) $. The relationship between the incident and the reflected geodesic arcs is then given by the transformation 
\begin{equation}
((\theta_i , \phi_i), (\theta_f , \phi_f)) \mapsto ((\theta'_i , \phi'_i), (\theta'_f , \phi'_f)) = T((\theta_i , \phi_i), (\theta_f , \phi_f)).
\end{equation}
We denote $T$ as the bounce map.

To make the action of the bounce map more explicit, consider a geodesic arc $((\theta_i , \phi_i), (\theta_f , \phi_f))$ that intercepts a face $F$ of $\Pi$ at a point $a$. As $F$ is a two-dimensional orientable surface, we can choose the inward unit normal $\eta_a$ to $F$ at $a$. We then select the unique geodesic that originates at $a$ and lies within the plane spanned by the incident ray $((\theta_i , \phi_i), (\theta_f , \phi_f))$ and $\eta_a$, which also represents the reflection of the incident ray in the hyperplane defined by $\eta_a$. The reflected ray $((\theta'_i , \phi'_i), (\theta'_f , \phi'_f))$ departs from the point $a$.\\

A billiard trajectory within $\Pi$ is thus a curve, parameterized by arc-length, consisting of the geodesic arcs reflected off the faces of the polyhedron. Formally, such a trajectory can be represented as a sequence $\gamma = ((\theta^n_i , \phi^n_i), (\theta^n_f , \phi^n_f))_{n \in \mathbb{Z}}$ where each element $((\theta^n_i , \phi^n_i), (\theta^n_f , \phi^n_f))$ is obtained by applying the bounce map $T$ to the preceding element: \begin{equation}
((\theta^n_i , \phi^n_i), (\theta^n_f , \phi^n_f)) = T ((\theta^{n-1}_i , \phi^{n-1}_i), (\theta^{n-1}_f , \phi^{n-1}_f)).
\end{equation} 
Note that we exclude billiard trajectories that intersect the vertices or edges of $\Pi$ from our consideration.

\subsection{Tessellating the Hyperbolic Space and the Unfolding of Billiard Trajectories}

Tessellations represent an essential tool in the analysis of hyperbolic spaces, both in the Poincaré ball model $\mathbb{B}^3$ and the hyperboloid model $\mathbb{H}^3$. We will next explore the concept of tessellation, which will enable us to unfold the billiard trajectories within a polyhedron.\\

A \emph{tessellation} of $\mathbb{B}^3$ is a subdivision of $\mathbb{B}^3$ into a collection of ideal polyhedron tiles $\Pi_i, \ i \in \Lambda$ following certain rules. Firstly, every point $a \in \mathbb{B}^3$ lies within some polyhedron, that is $\forall\ a \in \mathbb{B}^3,\ \exists\ i \in \Lambda $ with $a \in \Pi_i$. Secondly, the intersection of any two polyhedra is either empty, a single vertex, a single edge or an entire face. Formally, $\forall\ i \neq j,\ \Pi_i \cap \Pi_j$ is either $\emptyset$, a single vertex common to both, a single common edge or an entire common face. Lastly, for any two polyhedra, there exists an isometry $f_{i,j}$ of $\mathbb{B}^3$  such that $f_{i,j}(\Pi_i) = \Pi_j$. These rules apply analogously for a tessellation of $\mathbb{H}^3$.\\ 

Given an ideal polyhedron $\Pi$ in the hyperbolic space, we can produce a tessellation by reflecting $\Pi$ across each of its faces, and then repeating the process for the resulting reflections, indefinitely. This process yields a collection of ideal polyhedrons that cover $\mathbb{B}^3$, thereby tessellating it. This method is known as the \emph{Katok-Zemlyakov unfolding method}, developed by \texttt{A.B. Katok} and \texttt{A.N. Zemlyakov} \cite{katok}. In the case of a tetrahedron, its faces can be labeled arbitrarily as $1, 2, 3, 4$. When reflecting about a face $i$, the labels on the resulting faces are transformed to $1^i,2^i,3^i, 4^i$. This notation is retained consistently for subsequent reflections. 

The unfolding method transforms an ideal polyhedron in the hyperbolic space into a sequence of isometric polyhedrons, akin to a \say{tube.} Within this tube, a billiard trajectory appears as a full geodesic. This representation enables us to \emph{unfold} a billiard trajectory in $\Pi$ into an uninterrupted geodesic, providing an intuitive visualization of the trajectory's evolution.\\

Symbolic Dynamics provides a robust framework to analyze the long-term behaviour of billiard trajectories. This field involves discretizing the state space of a dynamical system, replacing continuous trajectories with sequences of symbols that represent states of the system at discrete time intervals. Despite losing some detailed information, this discretization facilitates the analysis by simplifying the trajectories.

To formalize this concept, we start with a finite set $\mathcal{A}$ of symbols, which we call the \emph{alphabet}. The \emph{full $\mathcal{A}$ shift} is the set of all bi-infinite sequences of symbols from $\mathcal{A}$, denoted by $ \mathcal{A}^{\mathbb{Z}} = \{ x = \ldots x_{-1}.x_0x_1\ldots  :\ x_i\ \in\ \mathcal{A}\ \forall\ i \in \mathbb{Z} \}$. We endow this set with the product topology, where the distance between two sequences $x = \ldots x_{-1}.x_0x_1\ldots  $ and $y = \ldots y_{-1}.y_0y_1\ldots  \in \mathcal{A}^{\mathbb{Z}}$ is defined by
\begin{equation}\label{metric}
d(x,y) = \inf \left\lbrace  \frac{1}{2^m}  : x_n = y_n \ \text{for} \ |n| < m  \right\rbrace.
\end{equation}

The dynamics on the full shift is defined by the \emph{shift map} $\sigma$, which shifts the elements of a sequence forward by one, i.e., $ (\sigma(x))_i\ =\ x_{i+1}$. A \emph{shift space} is then a subset $X \subseteq \mathcal{A} ^ {\mathbb{Z}}$ that is closed and invariant under the shift map. For more detailed exploration of Symbolic Dynamics, refer to \cite{marcus,morse}.

\section{Pointed Geodesics and Billiards in an Ideal polyhedron}

The concept of pointed geodesics was first introduced in \cite{nagar} in the context of billiards inside certain classes of polygons in the hyperbolic plane. Here, we extend these definitions to the hyperbolic space.

\subsection{Pointed Geodesics}

\bd Let $\gamma = ((\theta^n_i , \phi^n_i), (\theta^n_f , \phi^n_f)) \label{key}_{n \in \mathbb{Z}}$ be a billiard trajectory in a polyhedron $\Pi$ in $\mathbb{B}^3$. Then for a fixed $n \in \mathbb{Z}$, we call $((\theta^n_i , \phi^n_i), (\theta^n_f , \phi^n_f)) $ a base arc of $\gamma$.\ed 

Note that every base arc is a compact subset of $\mathbb{B}^3$ and in turn uniquely determines the billiard trajectory under the restrictions imposed by the reflection rules. 

\bigskip
\bd For a given base arc $((\theta_i , \phi_i), (\theta_f , \phi_f)) $ defining $\gamma $, we will call $(\gamma,((\theta_i , \phi_i),\\ (\theta_f , \phi_f))) $ a \emph{pointed geodesic}. \ed 

This implies that a pointed geodesic $(\gamma, ((\theta_i , \phi_i), (\theta_f , \phi_f)))$ is identified with 
\begin{equation}
\ldots  (T^{-1}((\theta_i , \phi_i), (\theta_f , \phi_f))).((\theta_i , \phi_i), (\theta_f , \phi_f))(T((\theta_i , \phi_i), (\theta_f , \phi_f))) \ldots  \in \K(\mathbb{B}^3)^{\Z}
\end{equation}
by fixing the position of the base arc $((\theta_i , \phi_i), (\theta_f , \phi_f)) $ in the bi-sequence. Here for a metric space $(X,d)$, we  denote
the space of all compact subsets of $X$ by  $\K(X)$, endowed with the Hausdorff topology. We label the faces of $\Pi$ with letters \emph{1,\ldots ,k} arbitrarily. A pointed geodesic can be encoded naturally by collecting the labels in the order in which it hits the faces of the polyhedron $\Pi$, by marking the face hit by the base arc and then reading the past and future hits. Thus every pointed geodesic produces a  bi-infinite sequence $\ldots a_{-1}.a_0a_1\ldots $ with $a_j\ \in\ \{1,\ldots ,k \}$.
 
 \bigskip

 \bd  \text{We call} $ \mathbf{G} = \mathbf{G}_{\Pi}= \{ \big(\gamma, ((\theta_i , \phi_i), (\theta_f , \phi_f)) \big) : \gamma  = \big( T^n((\theta_i , \phi_i), (\theta_f , \phi_f)) \big)_{n \in \mathbb{Z}} \} $, the space of all pointed geodesics on $\Pi$.
 
 Here, $ T^0 ((\theta_i , \phi_i), (\theta_f , \phi_f))$ simply denotes $((\theta_i , \phi_i), (\theta_f , \phi_f)) $. $\mathbf{G} \subset \K({\mathbb{B}^3})$ and so    ${\mathbf{G}}$ can be endowed with the natural Hausdorff metric $d_H$, and in turn inherits the Hausdorff topology. \ed

We define a function $ d_{\mathbf{G}} : \mathbf{G} \times \mathbf{G} \to \mathbb{R}$ by:
 \begin{equation} 
 \begin{aligned}
 \label{dg}
 d_{\mathbf{G}} \Big(\big(\gamma, ((\theta_i , \phi_i), (\theta_f , \phi_f)) \big),&\big(\gamma', ((\theta'_i , \phi'_i), (\theta'_f , \phi'_f)) \big)\Big) \\                      &= 
 \max\{ d_{\partial\mathbb{B}^3}((\theta_i,\phi_i ),(\theta'_i, \phi'_i)),d_{\partial\mathbb{B}^3}((\theta_f,\phi_f ),(\theta'_f, \phi'_f))\}. 
\end{aligned}
 \end{equation}

\bigskip

\bpr Suppose that $\mathbf{G}$ is the space of pointed geodesics associated with a polyhedron $\Pi$ in $\mathbb{B}^3$, then $d_{\mathbf{G}}$ is a metric on ${\mathbf{G}}$. \epr

\bo  By definition, $d_{\mathbf{G}}$ is non-negative.
 \begin{equation} 
d_{\mathbf{G}} \Big(\big (\gamma, ((\theta_i , \phi_i), (\theta_f , \phi_f)) \big),\big(\gamma', ((\theta'_i , \phi'_i), (\theta'_f , \phi'_f)) \big)\Big)\\ = 0,
 \end{equation}
implies 
 \begin{equation} 
 d_{\partial\mathbb{B}^3} ((\theta_i , \phi_i), (\theta'_i , \phi'_i))  = 0 , d_{\partial\mathbb{B}^3} ((\theta_f , \phi_f), (\theta'_f , \phi'_f)) = 0.
  \end{equation}

This means that 
\begin{equation} 
(\theta_i , \phi_i)= (\theta'_i , \phi'_i), (\theta_f , \phi_f)= (\theta'_f , \phi'_f),
 \end{equation}
which further gives 
\begin{equation} 
((\theta_i , \phi_i), (\theta_f , \phi_f)) = ((\theta'_i , \phi'_i), (\theta'_f , \phi'_f)).
\end{equation}
If the base arcs corresponding to two pointed geodesics are coincident, then their associated trajectories are also same due to the constraints of the reflection rules. Thus, 
\begin{equation}
\Big(\gamma, \big((\theta_i , \phi_i), (\theta_f , \phi_f)\big) \Big) =\Big(\gamma', \big((\theta'_i , \phi'_i), (\theta'_f , \phi'_f)\big) \Big).
\end{equation}
The  triangle inequality and symmetry for $d_{\mathbf{G}}$ follows from the respective properties of $d_{\partial\mathbb{B}^3}$ which establishes  $d_{\mathbf{G}}$ as a metric on ${\mathbf{G}}$.\eo 
Next, we will show that the Hausdorff topology on $\mathbf{G}$ is same as the topology on ${\mathbf{G}}$ given by $d_{\mathbf{G}}$. The metric $d_H$ on $\mathbf{G}$ can be given as follows:
\begin{equation} 
\begin{aligned}
d_H & \Big(\big (\gamma, ((\theta_i , \phi_i), (\theta_f , \phi_f)) \big),\big(\gamma', ((\theta'_i , \phi'_i), (\theta'_f , \phi'_f)) \big)\Big)\\ := & d_H  \Big( \big((\theta_i , \phi_i), (\theta_f , \phi_f) \big),\big((\theta'_i , \phi'_i), (\theta'_f , \phi'_f) \big)\Big) \\
=  & \max\lbrace \displaystyle \sup_{Q \in  \big((\theta_i , \phi_i), (\theta_f , \phi_f) \big)} d\big(Q, \big((\theta'_i , \phi'_i), (\theta'_f , \phi'_f) \big)\big), \displaystyle \sup_{Q \in  \big((\theta'_i , \phi'_i), (\theta'_f , \phi'_f) \big)} d\big(Q, \big((\theta_i , \phi_i), (\theta_f , \phi_f) \big)\big)\rbrace
\end{aligned}
 \end{equation}

This definition holds because $\gamma$ is described by its base arc uniquely. If $\gamma \neq \gamma'$ then 
\begin{equation}
d_H  \Big(\big (\gamma, ((\theta_i , \phi_i), (\theta_f , \phi_f)) \big),\big(\gamma', ((\theta'_i , \phi'_i), (\theta'_f , \phi'_f)) \big)\Big) >0.
\end{equation}
This notion of distance between two pointed geodesics is the Hausdorff distance between the corresponding base arcs. We will represent the space of all base arcs on a polyhedron $\Pi$ associated with billiard trajectories by $\mathbf{B}(\Pi)$ or $\mathbf{B}$. As was the case in $2$-dimensions, $\mathbf{B} \subset \K(\mathbb{B}^3)$, where $\K(X)$ denotes the space of all compact subsets of $X$, equipped with the Hausdorff topology.
Therefore, again in this case, the Vietoris topology and Hausdorff topology match on $\mathbf{B}$ giving us a natural isometry between $(\mathbf{G},d_H)$ and $(\mathbf{B},d_H)$ for any ideal polyhedron $\Pi$. This gives us a one-one correspondence between the Vietoris topology on $\mathbf{B}$ and the topology generated by $d_H$ on $\mathbf{G}$ leading to  $(\mathbf{G},d_{\mathbf{G}})$ and $(\mathbf{B},d_{\mathbf{G}})$ being isometric.\\

In the earlier discussion, we defined the metric $d_{\mathbf{G}}$ on the space of pointed geodesics, $\mathbf{G}$. Next we wish to show that the topology induced by $d_{\mathbf{G}}$ on $\mathbf{G}$ is equivalent to the topology inherited from the Hausdorff metric, $d_H$, on the space of compact subsets of $\mathbb{B}^3$. 

To that end, we need to demonstrate that for any pointed geodesic $\gamma \in \mathbf{G}$, an open ball in $(\mathbf{G}, d_{\mathbf{G}})$ containing $\gamma$ also contains an open ball in $(\mathbf{G}, d_H)$ containing $\gamma$, and vice versa. In other words, we want to show that any open set in one topology contains an open set in the other, implying the two topologies are equivalent.\\

First, let's define the Hausdorff distance on the space of compact subsets of $\mathbb{B}^3$ more rigorously. Let $A$ and $B$ be compact subsets of $\mathbb{B}^3$. For a set $A$ and a point $b \in \mathbb{B}^3$, define the distance from $b$ to $A$ as 
\begin{equation}
    d(b,A) = \inf_{a \in A} d(b,a).
\end{equation}

Then, the Hausdorff distance between $A$ and $B$ is given by
\begin{equation}
    d_H(A,B) = \max \left\lbrace \sup_{a \in A} d(a,B), \sup_{b \in B} d(b,A) \right\rbrace.
\end{equation}

In the context of pointed geodesics, each $\gamma \in \mathbf{G}$ can be viewed as a compact subset of $\mathbb{B}^3$ (since it is a geodesic segment in the compact ball $\mathbb{B}^3$), and thus we can apply $d_H$ to pairs of pointed geodesics.

\bigskip 

\bt Suppose $\mathbf{G}$ be the space of pointed geodesics on a polyhedron $\Pi$ in $\mathbb{B}^3$, then  $d_{\mathbf{G}}$ and $d_H$ are topologically consistent on $\mathbf{G}$. \et
\bo \begin{figure}[ht!]
\centering
\includegraphics[width=15cm,height=12.7cm]{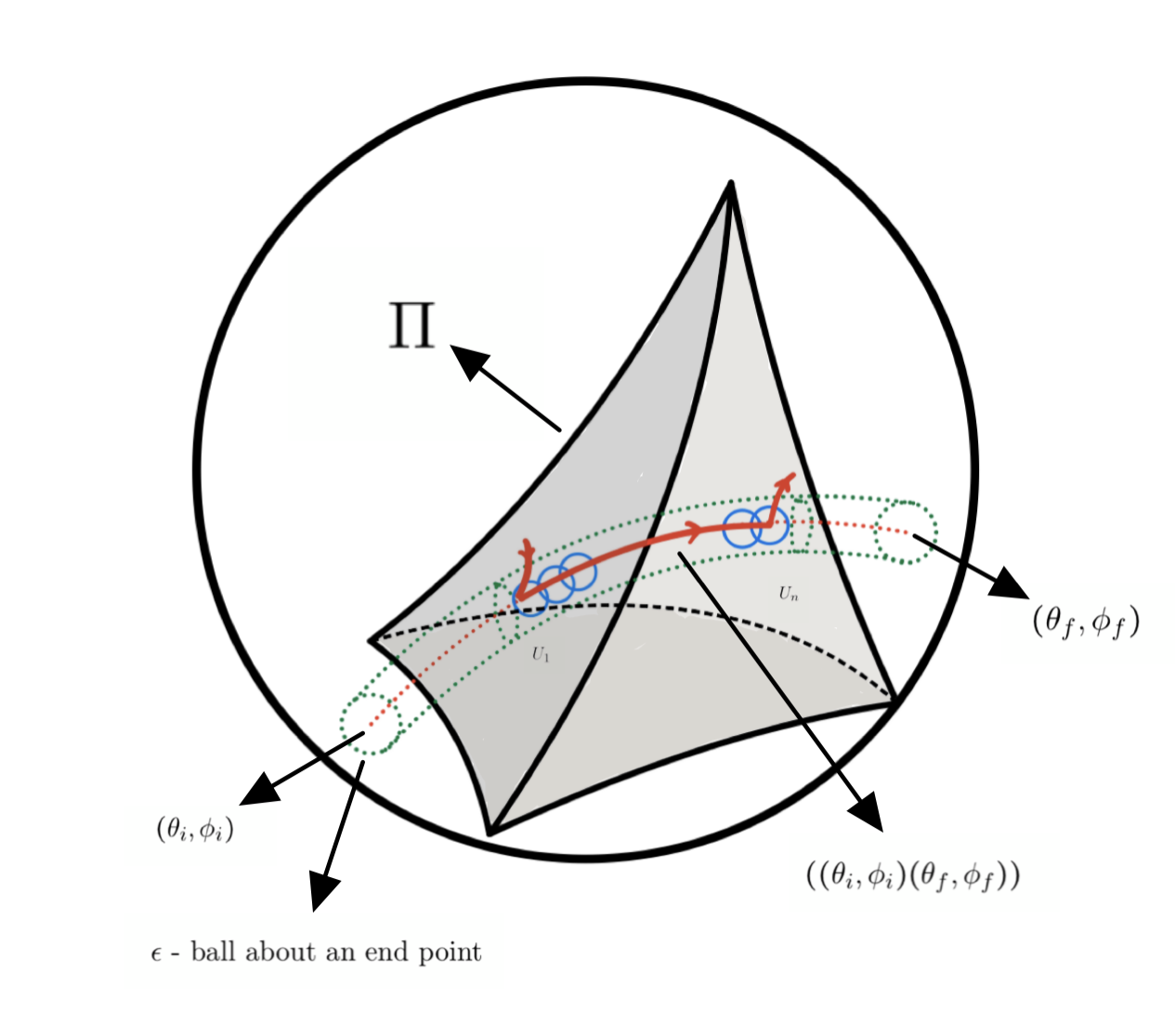}
\caption{An $\epsilon$-ball about a pointed geodesic}
\end{figure}
It is sufficient to prove that $d_{\mathbf{G}}$ and $d_H$ are topologically same on $\mathbf{B}$. The space $\mathbf{B}$ via the metric $d_H$ gets the induced topology given by $\mathbf{B}\ \subset\ \K({\mathbb{B}^3})$.
For $\epsilon > 0$, take
\begin{equation}
V = \big\{ ((\theta'_i, \phi'_i),(\theta'_f,\phi'_f)): d_{\mathbf{G}} \Big(((\theta'_i, \phi'_i),(\theta'_f,\phi'_f)), \big ((\theta_i, \phi_i),(\theta_f,\phi_f)) \Big) < \epsilon\big\},
\end{equation}
implying  
\begin{equation}
d_{\partial\mathbb{B}^3} ((\theta'_i, \phi'_i),\\(\theta_i,\phi_i)), d_{\partial\mathbb{D}} ((\theta'_f, \phi'_f),(\theta_f,\phi_f)) < \epsilon.
\end{equation}
Without any loss of generality, we can assume that $\epsilon$ is small enough such that the $\epsilon$-tube of the base arcs about $((\theta_i, \phi_i),(\theta_f,\phi_f)) $ doesn't contain any vertex or edge of $\Pi$, as has been shown in figure 2. Take the open balls $U_1,U_2,\ldots ,U_n$ in $\mathbb{B}^3$ such that 
\begin{equation}
((\theta_i, \phi_i),(\theta_f,\phi_f)) \subset \cup_{i=1}^{n} U_i,\ ((\theta_i, \phi_i),(\theta_f,\phi_f)) \cap U_i \neq \emptyset\ \forall\ i= 1,\ldots ,n 
\end{equation}
and each $U_i$ lying inside the $\epsilon$-tube. With this $<U_1,\ldots ,U_n>$ is open in $\K({\mathbb{B}^3})$, which means that $\mathbf{B}\ \cap <U_1,\ldots ,U_n>$ is open in $\mathbf{B}$ and is completely inside the $\epsilon$-tube. Thus 
\begin{equation}
((\theta_i, \phi_i),(\theta_f,\phi_f)) \in \mathbf{B}\ \cap<U_1,\ldots ,U_n>\ \subset V.
\end{equation}

Conversely, consider a basic open set $\mathbf{B}\ \cap <U_1,\ldots ,U_n>$ containing a base arc $((\theta_i, \phi_i),(\theta_f, \phi_f)) $.
Without any loss of generality, we can take $U_i's$ as the open discs in $\mathbb{B}^3$. Set 
\begin{equation}
W_{ij} = \{ p \in \mathbb{B}^3 : p \in U_i\ \cap\ U_j\  \forall i,j \in \{1,\ldots ,n\}, i \neq j \},
\end{equation}
where each $W_{ij}$ is either $\emptyset$ or has two points. Let 
\begin{equation}
W_0 = \{p \in \mathbb{B}^3 : p \in (U_i \cap(\partial \Pi)_k) \cup(U_i \cap(\partial \Pi)_{k+1}) \ \forall\ i = 1,\ldots ,n \}.
\end{equation}
$(\partial \Pi)_k$ denotes the face of the polyhedron with label $k$.  Set 
\begin{equation}
W= (\cup_{i,j=1,i \neq j}^{n} W_{ij}) \cup W_0 
\end{equation}
and pick $\delta < \inf_{p \in W} (d_{\partial \mathbb{B}^3}(p, ((\theta_i, \phi_i),(\theta_f,\phi_f)))).$
Therefore, the $\delta$-tube 
\begin{equation}
V = \{ ((\theta'_i, \phi'_i),(\theta'_f,\phi'_f)) : d_{\mathbf{G}}(((\theta'_i, \phi'_i),(\theta'_f,\phi'_f)), ((\theta_i, \phi_i),(\theta_f,\phi_f))) < \delta \}
\end{equation}
lies in $\mathbf{B}\ \cap <U_1,\ldots ,U_n>$ completely.


\eo

\subsection{Billiards in Ideal polyhedrons}

To give structure to our proof for the main result (Theorem \ref{coding}), we  split it into three parts, wherein we start with a lemma from elementary metric space theory. Then we follow it up with establishing of coding rules (Theorem \ref{ideal}) followed by defining the conjugacy between the space of codes and the space of pointed geodesics (Theorem \ref{coding}).

\bl \label{lemma} Let $x \in \mathbb{B}^3$ be a fixed point. Then every hyperbolic plane $p \in \mathbb{B}^3$ which does not contain $x$ divides $\mathbb{B}^3$ into two open half spaces $H^+_p$ and $H^-_p$ with $x \in H^-_p$. Further, let $(p_n)_{n\geq 0}$ be a sequence of hyperbolic planes with the additional properties that 
 $p_{n+1} \subset H^+_{p_n}$ for all $n$ and
 $d(p_n,x) \rightarrow \infty$.
Then the halfspaces $H^+_{p_n}$  determine a unique ideal point $\eta  \in \partial \mathbb{B}^3 $ and each geodesic originating from a point $y \in H^-_{p_0}$ and ending in $\eta$ penetrates successively once through each of the hyperbolic planes $(p_n)_{n\geq 0}$.\el

\bigskip

\bt \label{ideal}
For a fixed $k \in 2 \mathbb{N}+2$, let $\Pi \subset  \mathbb{B}^3$ be a $k$-faced ideal polyhedron with labelling as follows: Mark the faces $1,2,\ldots ,k$ in an arbitrary order and then the vertex defined by faces $i_1,\ldots ,i_m $ takes the label $i_1 \ldots  i_m$ for each $i_j  \in   \{1,\ldots ,k \}$. Let $\Omega_{ij}$ denote the interior angle between the pair of faces of $\Pi$ labelled $i $ and $j\ \forall\ i \neq j $ and $i,j \in \{1,\ldots ,k\}$. Further, assume that $\lambda_{ij} = \pi / \Omega_{ij} \in \mathbb{N}$ for each $i \neq j$ and $i,j \in \{1,\ldots ,k\}$. Then  an equivalence class $[\ldots a_{-1}.a_0a_1\ldots ]$ denoted $(a_j)$ with $\ldots a_{-1}.a_0a_1\ldots  \in \{1,\ldots ,k\}^{\mathbb{Z}}$ is in $S(\Pi)$, if and only if \\
(a) $a_j \neq a_{j+1} \forall j \in \mathbb{Z}$.\\
(b) $(a_j)$ does not contain more than $\lambda_{ij}$ repetitions of symbols $i$ and $j$ which are labels of two adjacent faces.\\
(c) $(a_j)$ does not contain an infinitely repeated sequence of labels of  faces meeting at a vertex.\\
Further, every equivalence class of such bi-infinite sequences corresponds to one and only one billiard trajectory.
\et

\bo
First we establish the necessity of (a), (b) and (c).
(a) follows from the fact that a billiard trajectory cannot hit the same hyperbolic plane twice because a geodesic cannot intersect a hyperbolic plane more than once. For (b), suppose $F_j,F_k$ be two adjacent faces of $\Pi$. Using a suitable isometry we can consider without loss of generality that these two faces are represented by the vertical hyperbolic planes $\{y=0\}$ and $\{ax+by=0\}$ in $\mathbb{H}^3$ for fixed $a,b \in \mathbb{R}$ with $a^2+b^2 \neq 0, a \neq 0$. Suppose (b) does not hold and a billiard trajectory not starting or ending in the vertex at $\infty$ of $\Pi$ hits faces $F_j,F_k$ more than $\lambda_{jk}$ times(say $\mu_{jk}$ times). Each segment of the corresponding billiard trajectory lying between $F_j$ and $F_k$ is part of a semicircle hitting $\{z=0\} $ plane orthogonally because it cannot be a straight line perpendicular to $\{z=0\}$ otherwise it will hit the vertex at $\infty$. On unfolding this part of the trajectory we get a part of an orthogonal semicircle. On projecting this part onto the $x-y$ plane we see that the projection is subtending an angle $\mu_{jk} \Omega_{jk} > \pi$ at origin which leads to a contradiction since the projection is a euclidean straight line in $x-y$ plane.\\
For (c), the following argument holds for the case where an arbitrary number of faces of $\Pi$ meet at a vertex. For illustration, let us consider in particular a situation where $F_j,F_k,F_l$ be three faces of $\Pi$ meeting at a vertex $v$. 
Using a suitable isometry we can consider without loss of generality that these three faces are represented by the vertical hyperbolic planes $\{y=0\}, \{ax+by=0\}$ and $\{cx+dy=0\}$ with the corresponding vertex at $\infty$. Let us label $\{y=0\}$ by $1$, $\{ax+by=0\}$ by $2$ and $\{cx+dy=0\}$ by $3$ without any loss of generality, as referred in figure 3. Suppose (c) does not hold i.e. there exists a billiard trajectory whose code contains an infinite word $w$ with $w_i \in \{1,2,3\}\ \forall \ i$. On unfolding the corresponding part of the billiard trajectory we get a part of the geodesic which is euclidean semicircle and orthogonal to $x-y$ plane. (b) ensures that the copies of $\Pi$ come out of any euclidean circle drawn on the plane $\{z=0 \}$. With this a geodesic in $\mathbb{H}^3$ can hit only finitely many copies of $\Pi$ generated while unfolding, which contradicts the infinite cardinality of $w$.\\

\begin{figure}[ht!]
\centering
\includegraphics[width=14.6cm,height=13.6cm]{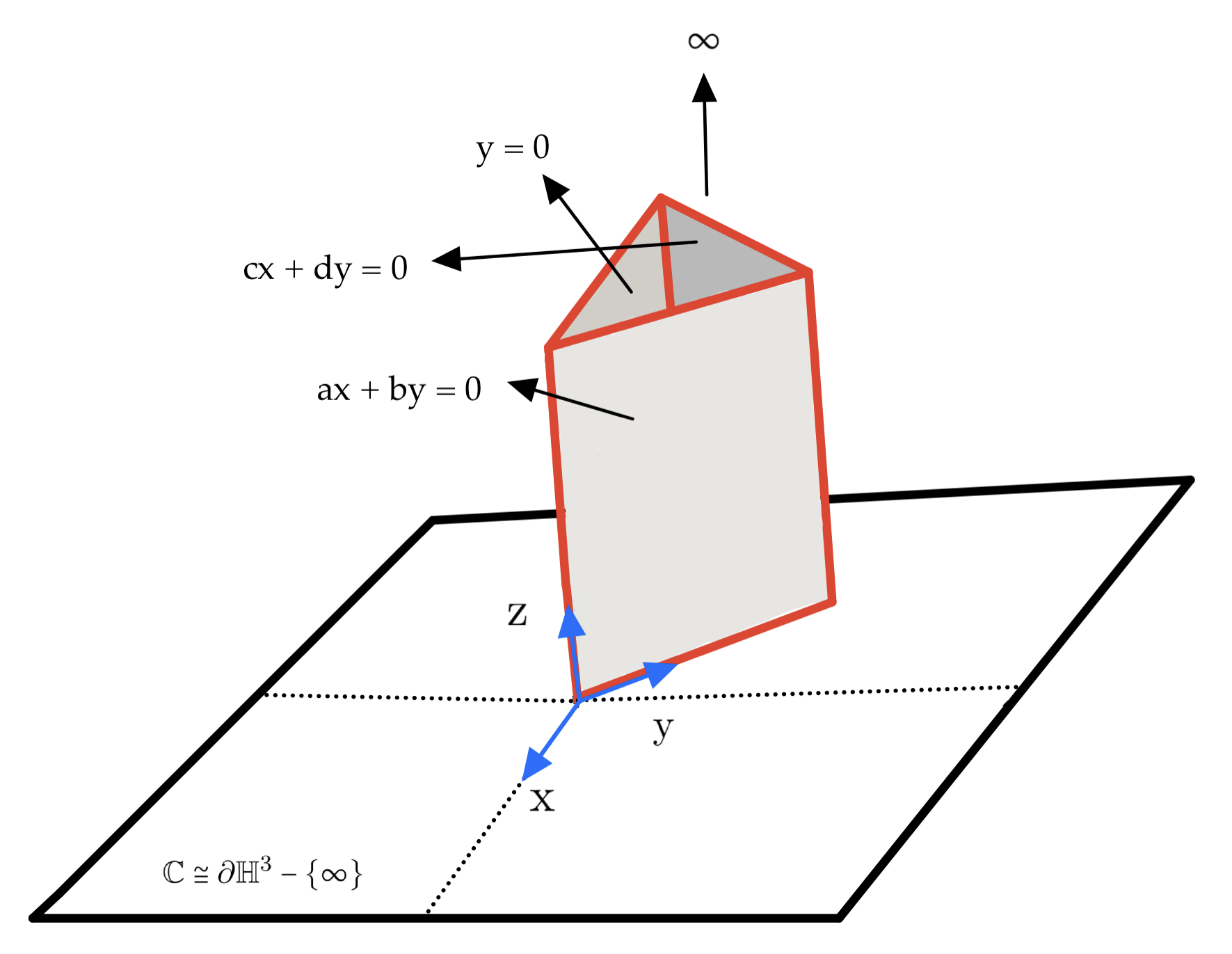}
\caption{An illustration of a vertex at $\infty$ defined by three vertically placed hyperbolic planes }
\end{figure} 

Conversely, we choose a sequence $(x_i)_{i \in \mathbb{Z}}$ which satisfies (a), (b) and (c) after fixing an ideal polyhedron $\Pi$ in $\mathbb{B}^3$ with faces labeled $1,2,3,4$ in an arbitrary order. Thus 
the chosen sequence $(x_i)_{i \in \mathbb{Z}}$ adheres to the restrictions imposed by (a), (b) and (c). We aim at creating an algorithm that starting from $(x_i)_{i \in \mathbb{Z}}$ generates the corresponding unique billiard trajectory in $\Pi$. We start by fixing an arbitrary point $A$ inside $\Pi$. The sequence $(x_i)_{i \in \mathbb{Z}}$ dictates us how to unfold the polyhedron $\Pi$ in $\mathbf{B}$. This procedure will produce copies of $\Pi = \Pi^0$ which we will label as $\Pi^{(i)}$ which serves as the reflected copy of $\Pi^{(i-1)}$ in the face labeled $x_i$ for $i \geq 1$ and the reflected copy of $\Pi^{(i+1)} $ in the face labeled $x_i$ for $i \leq -1$. Thus we get a bi-sequence $(\Pi^{(i)})_{i \in \mathbb{Z}}$ of isometric copies of $\Pi$. We will label the bi-sequence of faces in which the reflections are taking place as $(p_j)_{j \in \mathbb{Z}}$. Note that $p_j$ is labeled $x_j \forall j \in \mathbb{Z}$. The bi-sequence $(p_j)_{j \in \mathbb{Z}}$ obeys the condition (1) of the Lemma \ref{lemma}. Indeed, $A$ is not contained in any face $p_j$ and $p_{j+1} \subset H^+_{p_j}\ \forall\ j$.      Next, we establish condition (2) of Lemma \ref{lemma}. Let $\delta > 0$ be the least of all distances between all the non-adjacent planes of $\Pi$. Then $d(p_{j+1} , A) \geq d(p_{j+1}, p_j) + d(p_j, A) \geq \delta +d(p_j , A)$ for the non-adjacent planes $p_j , p_{j+1}$. Thus if we have infinitely many such pairs $(p_j , p_{j+1} )$ of non-adjacent planes for $j \geq 0$, we get $d(p_j , A) \rightarrow \infty$ which ensures the condition (2) of Lemma \ref{lemma}. This will define a unique limit point $\beta \in \partial \mathbb{B} $. The unique geodesic emanating from $\alpha$ and ending in $\beta$ is the unfolded trajectory corresponding to the code $(x_i)_{i \in \mathbb{Z}}$. If the non-adjacent pairs of planes are not infinite in either direction then we look for the triples $(p_{j-1}, p_j, p_{j+1})$  of non-adjacent planes of $\Pi$.

\eo

\bigskip
\bt \label{coding}
For $k \in 2\mathbb{N}+2$, let  $\Pi \subset \mathbb{B}^3$ be an ideal polyhedron with labelling as follows: Mark the faces $1,2,\ldots ,k$ in an arbitrary order and then the vertex defined by faces $i_1,\ldots ,i_m $ takes the label $i_1 \ldots  i_m$ for each $i_j  \in   \{1,\ldots ,k \}$. Let $\Omega_{ij}$ denote the interior angle between the pair of faces of $\Pi$ labelled $i $ and $j\ \forall\ i \neq j $ and $i,j \in \{1,\ldots ,k\}$. Further, assume that $\lambda_{ij} = \pi / \Omega_{ij} \in \mathbb{N}$ for each $i \neq j$ and $i,j \in \{1,\ldots ,k\}$.  Let $\mathbf{G}$ be the space of pointed geodesics on $\Pi$ and $X$ the space of all bi-infinite sequences $\ldots a_{-1}.a_0a_1\ldots  \in \{1,\ldots ,k\}^{\mathbb{Z}}$ satisfying (a), (b) and (c) from Theorem \ref{ideal}. Then $(\mathbf{G},\tau) \simeq (X,\sigma)$.
\et

\bo
Let us define $h : (\mathbf{G},\tau) \to (X,\sigma)$ by
\begin{equation}
\begin{aligned}
h  & \big(\gamma, (({\theta}_i,{\phi}_i),({\theta}_f , {\phi}_f)) \big) &= \ldots  a_{T^{-1} (({\theta}_i,{\phi}_i),({\theta}_f , {\phi}_f))}.a_{(({\theta}_i,{\phi}_i),({\theta}_f , {\phi}_f))} a_{T (({\theta}_i,{\phi}_i),({\theta}_f , {\phi}_f))}\ldots 
\end{aligned}
\end{equation}

Then
\begin{equation}
\begin{aligned}
 h & \big(\gamma, (({\theta}_i,{\phi}_i),({\theta}_f , {\phi}_f)) \big)  = h\big(\gamma', (({\theta'}_i,{\phi'}_i),({\theta'}_f , {\phi'}_f)) \big) \\
& \Rightarrow \ldots a_{T^{-1} (({\theta}_i,{\phi}_i),({\theta}_f , {\phi}_f))}.a_{(({\theta}_i,{\phi}_i),({\theta}_f , {\phi}_f))} a_{T (({\theta}_i,{\phi}_i),({\theta}_f , {\phi}_f))}\ldots  \\ &= \ldots a_{T^{-1} (({\theta'}_i,{\phi'}_i),({\theta'}_f , {\phi'}_f))}.a_{ (({\theta'}_i,{\phi'}_i),({\theta'}_f , {\phi'}_f))} a_{T (({\theta'}_i,{\phi'}_i),({\theta'}_f , {\phi'}_f))} \ldots  \\
& \Rightarrow \big(a_{T^n  (({\theta}_i,{\phi}_i),({\theta}_f , {\phi}_f))}\big)_{n \in \mathbb{Z}} = \big(a_{T^n  (({\theta'}_i,{\phi'}_i),({\theta'}_f , {\phi'}_f))}\big)_{n \in \mathbb{Z}}
\end{aligned}
\end{equation}

From Theorem \ref{ideal}, we have $ (T^n (({\theta}_i,{\phi}_i),({\theta}_f , {\phi}_f))\big)_{n \in \mathbb{Z}} = (T^n (({\theta'}_i,{\phi'}_i),({\theta'}_f , {\phi'}_f))\big)_{n \in \mathbb{Z}}$ 
and $\\ a_{ (({\theta}_i,{\phi}_i),({\theta}_f , {\phi}_f))} = a_{ (({\theta'}_i,{\phi'}_i),({\theta'}_f , {\phi'}_f))}$, which implies 
\begin{equation}
\big(\gamma,  (({\theta}_i,{\phi}_i),({\theta}_f , {\phi}_f)) \big) = \big(\gamma',  (({\theta'}_i,{\phi'}_i),({\theta'}_f , {\phi'}_f))\big).
\end{equation}
This gives the injectivity of h.
We get the surjectivity of $h$  again by using Theorem \ref{ideal} as each $(a_j)_{j \in \mathbb{Z}} \in S(\Pi)$ defines a unique billiard trajectory $\gamma$, which in turn implies that with the corresponding $\ldots a_{-1}.a_0a_1\ldots $, we get a unique base symbol $a_0$, which defines a base arc $ (({\theta}_i,{\phi}_i),({\theta}_f , {\phi}_f)) $ on $\gamma$, giving a unique pointed geodesic in $\mathbf{G}$. We have  
\begin{equation}
    h  \big(\gamma,  (({\theta}_i,{\phi}_i),({\theta}_f , {\phi}_f)) \big)\\ = \ldots a_{-1}.a_0a_1\ldots 
\end{equation}
Then
\begin{equation}
\begin{aligned}
 h & \circ \tau\Big(\big(\gamma,  (({\theta}_i,{\phi}_i),({\theta}_f , {\phi}_f)) \big)\Big)\\ &= h\bigg(\tau\Big(\big(\gamma,  (({\theta}_i,{\phi}_i),({\theta}_f , {\phi}_f)) \big)\Big)\bigg) \\  &= h\Big(\big(\gamma, T( ({\theta}_i,{\phi}_i),({\theta}_f , {\phi}_f)) \big)\Big)\\ &= h\Big(\big(\gamma,  (({\theta}_{1i},{\phi}_{1i}),({\theta}_{1f} , {\phi}_{1f})) \big)\Big)\\
&= \ldots a_{T^{-1} (({\theta}_{1i},{\phi}_{1i}),({\theta}_{1f} , {\phi}_{1f}))}.a_{(({\theta}_{1i},{\phi}_{1i}),({\theta}_{1f} , {\phi}_{1f}))} a_{T (({\theta}_{1i},{\phi}_{1i}),({\theta}_{1f} , {\phi}_{1f}))}\ldots \\  &=\ldots a_{T^{-1}T (({\theta}_i,{\phi}_i),({\theta}_f , {\phi}_f))}.a_{T (({\theta}_i,{\phi}_i),({\theta}_f , {\phi}_f))} a_{TT (({\theta}_i,{\phi}_i),({\theta}_f , {\phi}_f))}\ldots \\ &= \ldots a_{(({\theta}_i,{\phi}_i),({\theta}_f , {\phi}_f))}.a_{T (({\theta}_i,{\phi}_i),({\theta}_f , {\phi}_f))} a_{T^2 (({\theta}_i,{\phi}_i),({\theta}_f , {\phi}_f))}\ldots \\ &= \sigma \Big(h\big(\gamma, (({\theta}_i,{\phi}_i),({\theta}_f , {\phi}_f))\big)\Big)\\ &= \sigma \circ h \big(\gamma, (({\theta}_i,{\phi}_i),({\theta}_f , {\phi}_f)) \big) 
\end{aligned}
\end{equation}

 (Here  $(({\theta}_{1i},{\phi}_{1i}),({\theta}_{1f} , {\phi}_{1f}))$ is the reflected geodesic arc corresponding to the incident  ray  $(({\theta}_i,{\phi}_i),({\theta}_f , {\phi}_f)))$

$\implies\ h \circ \tau = \sigma \circ h$, which implies $h$ is a homomorphism.\\

Let $U = [x_{-m}\ldots x_{-1}.x_0\ldots x_m]$ in $(X, \sigma)$ be an open set in $\mathbb{G}$. For a pointed bi-sequence $x \in$ U, we have corresponding $(x_n)_{n \in \mathbb{Z}}$, which generates a billiard trajectory $\gamma$ using Theorem \ref{ideal}. By pointing out the base arc $ (({\theta}_i,{\phi}_i),({\theta}_f , {\phi}_f)) $ corresponding to symbol $x_0$, we get a pointed geodesic $(\gamma, (({\theta}_i,{\phi}_i),({\theta}_f , {\phi}_f)))$. Let's label its pointed bi-sequence by  $y=\ldots y_{-1}.y_0 y_1\ldots $. Since $x$ and $y$ belong to same equivalence class, there exists an $s$ such that 
\begin{equation}
y_{[s-m,s+m]} = x_{-m}\ldots x_{-1}.x_0\ldots x_m.
\end{equation}
Therefore, $(\gamma,T^{-s} (({\theta}_i,{\phi}_i),({\theta}_f , {\phi}_f)))$ has its associated pointed billiard bi-sequence 
\begin{equation*}
h(\gamma,T^{-s} (({\theta}_i,{\phi}_i),({\theta}_f , {\phi}_f)))\   \in U.
\end{equation*}
We will construct $m$ future and past copies of $\Pi$ in $\mathbb{B}^3$ by reflecting $\Pi$ in its faces under the order given by $h(\gamma,T^{-s} (({\theta}_i,{\phi}_i),({\theta}_f , {\phi}_f)))\ \in U.$ Label $T^{-s} (({\theta}_i,{\phi}_i),({\theta}_f , {\phi}_f)) $ as $ (({\theta'}_i,{\phi'}_i),({\theta'}_f , {\phi'}_f)) $.

Define $\delta_1$ as  
\begin{equation}
\delta_1 = \displaystyle\min_{i \in \{1,\ldots ,k\}}\big\{d_{\partial\mathbb{B}^3}(A^m_i,(\theta'_f,\phi'_f)),d_{\partial\mathbb{B}^3}(A^{-m}_i,(\theta_i', \phi_i'))\big\}. 
\end{equation}
    Choose $\epsilon$ such that $0 < \epsilon <\delta_1$. If  
\begin{equation}    
\big(\gamma', (({\theta'}_i,{\phi'}_i),({\theta'}_f , {\phi'}_f))\big) \in B_{\epsilon}\big(\gamma, (({\theta}_i,{\phi}_i),({\theta}_f , {\phi}_f))\big),
\end{equation}
    then 
  \begin{equation}    
     [h\big(\gamma', (({\theta'}_i,{\phi'}_i),({\theta'}_f , {\phi'}_f))\big)]_{[-m,m]} = \displaystyle x_{-m}\ldots x_{-1}x_0\ldots x_m. 
     \end{equation}
    Thus  
     \begin{equation} 
    h\big(\gamma', (({\theta'}_i,{\phi'}_i),({\theta'}_f , {\phi'}_f))\big) \in U
    \end{equation}
    and this implies that 
     \begin{equation}
    h\Big(B_{\epsilon}\big(\gamma, (({\theta}_i,{\phi}_i),({\theta}_f , {\phi}_f))\big)\Big)   \subseteq U.
    \end{equation}
    Therefore, $h$ is continuous.\\

Conversely, let 
 \begin{equation}
V = B_{\epsilon}\big(\gamma, (({\theta}_i,{\phi}_i),({\theta}_f , {\phi}_f))\big)
 \end{equation}
be open in $\mathbb{G}$. Thus $\big(\gamma', (({\theta'}_i,{\phi'}_i),({\theta'}_f , {\phi'}_f))\big)  \in V $ if and only if 
 \begin{equation}
d_{\partial\mathbb{D}}((\theta_i,\phi_i),(\theta'_i , \phi'_i)),d_{\partial\mathbb{D}} ((\theta_f,\phi_f),(\theta'_f , \phi'_f)) < \epsilon.
 \end{equation}
We can tesselate $\mathbb{B}^3$ with $\Pi$ and its copies generated by reflecting $\Pi$ about its sides and doing the same for the reflected copies along the unfolded geodesic generated by $\gamma$. Label the vertices of $\Pi$ arbitrarily by $A_1, A_2,\ldots ., A_k$ and the vertices of the $i^{th}$ copy of $\Pi$ by $A^i_1, A^i_2,\ldots , A^i_k$. Define $p$ to be the largest positive integer such that  $A^i_1, A^i_2,\ldots , A^i_k $ are not in $\epsilon$-ball about  $(({\theta}_i,{\phi}_i),({\theta}_f , {\phi}_f))$  for $i = -p,-p+1,\ldots ,0,1,\ldots ,p.$ which means that
 \begin{equation}
h^{-1}([x_{-p}\ldots x_{-1}x_0\ldots  x_k]) \subseteq V
 \end{equation}
and thereby $h^{-1}$ is continuous.
\eo

It is noted here that the space \emph{X} described above is not closed as its  limit points of type $w\overline{abc}$ and $\overline{abc} w$ where $a,b,c$ are symbols appearing on faces meeting at a common vertex, do not lie in \emph{X}. 
\begin{equation}
\begin{aligned}
X &= \{ \ldots x_{-1}.x_0 x_1\ldots   \in \{1,\ldots ,k\}^{\mathbb{Z}} : x_i \neq x_{i+1}\ \forall\ i\ \text{and}\ \ldots x_{-1}.x_0 x_1\ldots  \neq \ w\overline{abc} ,\ \overline{abc} w \\ & \text{for}\ \text{any}\  a,\ b ,\ c\ \in \{1,\ldots ,k\}\  \text{sharing the same vertex and word}\ w, \ldots x_{-1}.x_0x_1\ldots\\ &  \neq w(ab)^{\mu_{ab}}w',  \mu_{ab} > \lambda_{ab} \forall\ \text{labels}\
	a,b \  \text{sharing an edge   and arbitrary}\ w, w' \}
\end{aligned}	
\end{equation}

Therefore, we go further and define the closure of \emph{X} in $\{1,\ldots ,k\}^{\mathbb{Z}}$, labelling it $\tilde{X}$. We can split  $\tilde{X}$ as $X \cup X'$, where $X'$ is the set of all limit points of \emph{X}.
\begin{equation}
\begin{aligned}
	\tilde{X} = \{\ldots x_{-1}.x_0 x_1\ldots  \in \{1,\ldots ,k\}^{\mathbb{Z}}: x_i \neq x_{i+1} \forall i,
	 \ldots x_{-1}.x_0x_1\ldots  \neq w(ab)^{\mu_{ab}}w', \mu_{ab} > \lambda_{ab} \\ \forall\ \text{labels}\
a,b \ \text{sharing an edge   and arbitrary}\ w, w' \}
\end{aligned}
\end{equation}

Thereby, $\tilde{X}$ has a finite forbidden set and thus is an \emph{SFT}. This places \emph{X} densely inside $\tilde{X}$. 
$\tilde{X}$ being the completion of \emph{X} is also the compactification of \emph{X}. 

\section{Conclusion}

In this paper, we have presented a detailed study of billiards in ideal polyhedrons in the hyperbolic space $\mathbb{B}^3$. Through a novel and well-defined coding system for billiard trajectories, we were able to create a robust connection between the physicality of these trajectories and the abstract realm of symbolic dynamics. In our analysis, we've shown a topological equivalence between the space of pointed geodesics on the polyhedron and the space of bi-infinite sequences that satisfy certain constraints, providing a solid foundation for understanding the complexity of our dynamical system.

Notably, our exploration unveiled the fact that the space of our defined sequences, although not closed, can be densely embedded into a symbolic shift space which possesses beneficial properties of closure, compactness, and finiteness. This is of significant importance as it has allowed us to rigorously analyze the dynamics of our system using the framework of symbolic dynamics, a tool frequently used to dissect complex dynamical systems.\\

One might question the practicality of such an abstract construction; yet, the bridge we've constructed between the concrete and the abstract proves itself to be rather powerful. The interplay between the geometry of the polyhedron and the symbolic dynamics of the shift map on the space of sequences illuminates intricate details about the dynamics of billiard trajectories in a way that would be challenging with geometric or physical analysis alone.

A number of interesting prospects have emerged from our work. The properties of the symbolic shift map on the shift space, including its possible mixing or transitive properties and the existence of periodic points, present enticing areas for further exploration. These properties could yield valuable insights into the behaviour of billiard trajectories on the polyhedron. Furthermore, we are intrigued by the potential to delve deeper into the metric entropy of the shift map on our space, which could offer a quantitative measure of the \say{chaos} in our billiard system. Drawing from the forbidden sequences that we've established, the calculation of entropy could further extend our understanding of the complexity of this system. We also envisage applying advanced topics such as the Patterson-Sullivan measures associated with our dynamical system, and the encompassing study of the thermodynamic formalism of our system. By correlating the geometry of the hyperbolic space, the symbolic dynamics, and the Patterson-Sullivan measures, we can hope to build a comprehensive view of our system. The emergence of an underlying hyperbolic dynamical system from our symbolic shift map is another captivating prospect. The concept of stable and unstable manifolds at each point could bring a whole new perspective to our understanding of the system's behaviour over time.\\

In conclusion, while our study offers a rich and rigorous analysis of billiards in ideal polyhedrons in $\mathbb{B}^3$, it also uncovers a plethora of intriguing avenues for future exploration. The marriage of physical geometry and abstract symbolic dynamics has proven to be a potent tool for understanding such complex systems. We anticipate that the further application and exploration of the concepts and methods presented here will continue to unravel the intricate and fascinating dynamics of billiard trajectories in hyperbolic space.

\section*{Data Availability Statement}

Data sharing is not applicable to this article as no new data was created or analysed in this study.

\section*{Acknowledgements}
This work was completed with the support from CSIR, India under the SPM fellowship. The author wishes to   thank \textsc{Prof. Anima Nagar} for the discussions and providing many useful insights into the problem.

\end{document}